\documentclass[twoside, 11pt]{article}
\usepackage{amsmath} 
\usepackage[e]{esvect}
\usepackage{graphicx}
\usepackage{epstopdf}
\usepackage[adobe-utopia]{mathdesign}
\usepackage[T1]{fontenc}
\usepackage[english]{babel}
\usepackage{color}

\usepackage{hyperref}

\title{\textsc{No more than $2^{d+1}-2$ nearly neighbourly simplices in $\mathbb R^d$}}
\date{}
\author{{Andrzej P. Kisielewicz \& Krzysztof Przes\l awski}\\
{\small \sf  A.Kisielewicz@wmie.uz.zgora.pl \,\, K.Przeslawski@wmie.uz.zgora.pl}\\
{\small\it Wydzia{\l} Matematyki, Informatyki i Ekonometrii, Uniwersytet Zielonog\'orski,}\\
{\small\it ul. Podg\'orna 50, 65-246 Zielona G\'ora, Poland}}
\newtheorem{lemat}{\sc Lemma}
\newtheorem{tw}{\sc Theorem}

\newtheorem{wn}{\sc Corollary}
\newtheorem{df}{\sc Definition}

\newtheorem{uw}{\sc Remark}
\newtheorem{uwi}[uw]{\sc Remarks}

\newtheorem{nap}{\sc Example}

\newtheorem{nps}[nap]{\sc Examples}

\def\ka #1{\mathscr{#1}}
\def\kal #1 #2{\mathscr{#1}^{#2}}
\def\proof{\noindent \textit{Proof.\,\,\,}}

\def\er{\mathbb{R}}

\def\Aut #1 #2{\operatorname{Aut}^{#1} (#2)}

\def\skalar #1#2{\langle #1,#2\rangle}

\def\prop #1{\operatorname{prop}\,(#1)}
\def\jed #1{\boldsymbol{1}_{#1}}

\def\bred #1 {\colorbox{red}{ #1}}
\def\red #1 {{\color{red} #1 }} 

\graphicspath{{figpdfkolor/}}
\DeclareGraphicsExtensions{.pdf,.png} 

\begin{document}
\maketitle

\begin{abstract}
We prove a combinatorial theorem on families of disjoint sub-boxes of a discrete cube, which implies that there are at most $2^{d+1}-2$ nearly neighbourly simplices in $\er^d$. 
\end{abstract}

\section{Introduction}
A family of $d$-simplices in $\er^d$ is \textit{nearly neighbourly} if every two members are separated by a hyperplane that contains a facet of each. This notion is related to the more restrictive neighbourliness. Let us recall that a family of $d$-simplices is called \textit{neighbourly} if the common part of every two members is a $(d-1)$-dimensional set. 

It has been repeatedly conjectured that the maximum cardinality $c_d$ of a neighbourly family of $d$-simplices is $2^{d}$ (see \cite{Zaks1} for further references). The conjecture is verified only up to dimension $3$. F. Bagemihl \cite{Bagemihl} proved that $8\le c_3\le 17$. V. Baston \cite{Baston} proved $c_3\le 9$. The final step, $c_3=8$, was made by J. Zaks \cite{Zaks2}. The same author \cite{Zaks1} showed by a clever construction that $c_d\ge 2^d$. It was M. A. Perles \cite{Perles} who had verified $c_d\le 2^{d+1}$. A slightly better estimate $c_d\le 2^{d+1}-1$ was shown in \cite[Chapter 14]{AZ}. (This chapter, together with a recent post \cite{Kalai} on G. Kalai's blog, is a great introduction to the subject of (nearly) neighbourly families of simplices.) Recently, the present authors \cite{KP} improved this estimate by 1. This small progress has been achieved by using methods of discrete geometry. 

Even less is known about nearly neighbourly families of $d$-simplices. By the same argument as in \cite{AZ}, one can prove that there are no more than $2^{d+1}-1$ simplices in such a family. On the other hand, our reasoning in \cite{KP} does not extend to this case as it heavily depends on the fact that for every two members of a neighbourly family of simplices there is only one hyperplane separating them that contains a facet of each. The aim of the present note is to show that $2^{d+1}-2$ is also an upper bound for the cardinality of a nearly neighbourly family of $d$-simplices. This result seems to be new, excluding the case $d=3$, which was resolved by S. Furino, B. Gamble and J. Zaks (see \cite{FGZ} or \cite{Zaks2}). Their long and involved proof utilises the technique developed by H. Tverberg (see \cite{Tverberg} or \cite{AZ}) to prove the Graham-Pollak theorem. Our proof is considerably shorter, and works in all dimensions.

\section{Main result}

Let $B=B_1\times\cdots \times B_n$ be a box contained in $\{0,1\}^n$. We define the following subset of $[n]=\{1,2,\ldots, n\}$:
$$
\prop B=\{i \colon B_i\neq \{0,1\}\}.
$$
\begin{tw}
\label{podstawa}
Let $3\le k < n$ be integers. Let $\ka B$ be a family of boxes in $\{0,1\}^n$ satisfying the following conditions:
\begin{enumerate}
\item[$(\alpha)$]
$|\prop B|=k$ for $B \in \ka B$;
\item[$(\beta)$]
$\prop A\neq \prop B$, whenever $A, B\in \ka B$ and $A\neq B$;
\item[$(\gamma)$] 
boxes belonging to $\ka B$ are pairwise disjoint.
\end{enumerate}
Then $|\ka B|\le 2^k-2$.
\end{tw}

Let us note that $(\gamma)$ means there is $i$ such that $A_i=\{0\}$ and $B_i=\{1\}$, or vice versa.

Let us recall some basic facts concerning Fourier transforms before going into the proof. Further details can be found in \cite{O'Donnell} or \cite{Terras}.

Let $V$ be the Euclidean vector space of all real-valued functions $f\colon\{0,1\}^n\to \er$ with scalar product defined by
$$
\skalar f g =\frac 1{2^n} \sum_{x\in \{0,1\}^n} f(x)g(x).
$$

For every $S\subseteq [n]$, let us define $\chi_S\colon \{0,1\}^n\to \er$ by 
$$
\chi_S(x)=(-1)^{\sum_{i\in S} x_i}.
$$ 
The functions $\chi_S$ form an orthonormal basis in $V$. Let $\hat{f}(S)$ be the coefficients of the Fourier expansion of $f\colon\{0,1\}^n\to \er $ with respect to this basis; that is,
$$
\hat{f}(S)=\skalar f{\chi_S}.
$$
This defines the Fourier transform $\hat f\colon 2^{[n]}\to \er$ of $f$. 

For every pair of functions $F,G\colon 2^{[n]}\to \er $, we define their convolution
$$
F*G(S)=\sum_{T\in 2^{[n]}} F(S\triangle T)G(T),
$$
where $\triangle$ denotes the symmetric difference. Recall that if $fg$ is the pointwise product of $f,g\in V$, then $\widehat{fg}=\hat f*\hat g$ (see e.g. \cite[Theorem 1.27]{O'Donnell}). In particular, if $f$ is an indicator function (that is, $f^2=f$), then 
\begin{equation}
\label{idempotent}
\hat f=\hat f*\hat f.
\end{equation}
We shall also need the following lemma:
\begin{lemat} 
\label{grupa}
Let $k$ be a positive integer and let $v$ be the $2$-adic order of $k$. If $\ka G\subseteq 2^{[n]}$ is a group under the operation of symmetric difference and $|C|=k$ for every $C\in \ka G'=\ka G\setminus \{\emptyset\}$, then 
\begin{equation}
\label{nier}
|\ka G|\le 2^{v+1}.
\end{equation}
\end{lemat}
\proof 
We may suppose that $|\ka G|\ge 4$, otherwise there is nothing to prove. Let us fix $K\in \ka G'$. Since there is yet another element $C$ in $\ka G'$, we deduce that $K\triangle C$ belongs to $\ka G'$. Since all these elements are of the same cardinality $k$, it follows that $k$ is even. If we interpret $\ka G$ as a vector space over $\mathbb F_2$, then it is easily seen that there is a subgroup $\ka H$ of $\ka G$ such that $\ka G$ decomposes into the direct sum
$$
\ka G=\ka H\oplus \{\emptyset , K\}. 
$$
Let us consider the mapping $\ka H\ni L\stackrel{\varphi}{\mapsto} L\cap K$. Observe that it is a group isomorphism onto $\varphi(\ka H)$. Indeed, it is a homomorhism, and if $L\neq \emptyset$, then $|L\cap K|=\frac k 2$, which implies injectivity of $\varphi$. Let $\varphi(\ka H)'=\varphi(\ka H)\setminus \{\emptyset\}$. Since every element of $\varphi(\ka H)'$ is of cardinality $\frac k 2$, we deduce by induction that 
$$
|\ka G|=2|\varphi(\ka H)|\le 2\times 2^v=2^{v+1}.
$$
\hfill$\square$

\medskip
\noindent
\textit{Proof of the theorem.}
Let $f=\sum_{B\in \ka B}\jed B$. By $(\alpha)$ and $(\beta)$,
$$
\hat{f}(\prop B)\in \left\{-\frac 1{2^k}, \frac 1{2^k}\right\},
$$ 
for every $B\in \ka B$. 
Let $m=|\ka B|$. Then
$$
\hat{f}(\emptyset)=\frac{m}{2^k}.
$$
By $(\gamma)$ and Bessel's inequality,
$$
\frac m {2^k}=\skalar f f \ge \hat{f}(\emptyset)^2+\sum_{B\in \ka B} \hat{f}(\prop B)^2= \frac {m^2+m}{2^{2k}}.
$$
Therefore, 
$$
m\le 2^k-1 
$$
and equality holds if and only if the support of $\hat f$ is equal to $\ka M=\{\emptyset\}\cup \{\prop B\colon B\in \ka B\}$. Let us fix $C\in \ka B$, and set $K=\prop C$. Let us suppose that $\ka M$ is indeed the support of $\hat f$. Since $(\gamma)$ implies that $f$ is an indicator function, we get by (\ref{idempotent})
$$
\hat{f}(K)=\sum_{Y\in \ka M} \hat{f}(K\triangle Y)\hat{f}(Y).
$$

Let $\ka N=\ka M\setminus\{\emptyset, K\}$. Then the preceding equation can be rearranged as follows
$$
1 =  2\hat f(\emptyset) +\frac 1{\hat{f}(K)}\sum_{Y\in \ka N} \hat{f}(K\triangle Y)\hat{f}(Y).
$$
Consequently, 
$$
1-\frac{1}{2^{k-1}}=2\hat f(\emptyset)-1=-\frac 1{\hat{f}(K)}\sum_{Y\in \ka N} \hat{f}(K\triangle Y)\hat{f}(Y).
$$
Observe that $\hat{f}(K\triangle Y)\in \{-\frac 1{2^k}, 0,\frac{1}{2^k}\}$ for every $Y\in \ka N$. Therefore,
$$
1-\frac{1}{2^{k-1}}\le \frac 1{|\hat{f}(K)|}\sum_{Y\in\ka N} |\hat{f}(K\triangle Y)\hat{f}(Y)|\le \frac{m-1}{2^k}= 1-\frac{1}{2^{k-1}}.
$$
This readily implies that $K \triangle Y$ is a nonempty element of the support of $\hat f$ for every $Y\in \ka N$. Thus, we have proved that $\ka M$ is a group under the symmetric difference. This group satisfies the assumptions of Lemma \ref{grupa} and consequently we deduce that if $v$ is a 2-adic order of $k$, then
$$
|\ka M|\le 2^{v+1} \le 2k< 2^k,
$$
which is a contradiction. 
\hfill $\square$

\medskip
Now, we demonstrate how to use Theorem \ref{podstawa} to show the corresponding result for nearly neighbourly simplices. We shall follow the approach that originates from \cite{Baston} and further developed in \cite{Perles}.

Let $\ka F$ be a nearly neighbourly family of $d$-simplices in $\er^d$, $d\ge 2$. Let us arrange all the hyperplanes spanned by the facets of simplices belonging to $\ka F$ into a sequence $H_1, \ldots, H_n$. Each $H_i$ splits $\er^d$ into two halfspaces. Let us call them $H_i^0, H_i^1$. For every simplex $\sigma\in \ka F$, let us define a unique box $B=B(\sigma)=B_1\times B_2\times\cdots\times B_n\subseteq \{0,1\}^n$ as follows
$$
B_i=
\begin{cases}
\{0\}, & \text{if $H_i$ is spanned by a facet of $\sigma$ and $\sigma\subset H_i^0$};\\
\{1\}, & \text{if $H_i$ is spanned by a facet of $\sigma$ and $\sigma\subset H_i^1$};\\
\{0,1\}, & \text{otherwise.}
\end{cases}
$$
Let $\ka B=\{B(\sigma)\colon \sigma \in \ka F\}$. It is clear that the family of boxes $\ka B$ satisfies the assumptions of Theorem \ref{podstawa} with $k=d+1$. Therefore, $|\ka B|\le 2^{d+1}-2$ and we have:
\begin{wn}
For every nearly neighbourly family $\ka F$ of $d$-simplices in $\er^d$
$$
|\ka F|\le 2^{d+1}-2.
$$
\end{wn}
(This holds true for $d=1$ as well.) 
\section{Remarks} 
\textbf{1.} For a fixed pair $(k, n)$, $3\le k<n$, let us consider all families $\ka B$ of boxes satisfying the assumptions of Theorem \ref{podstawa}. Let $b_{k,n}$ be the maximum of the numbers $ |\ka B|$. Let $b_k=\max_{n}b_{k,n}$. We show by means of examples that the upper bound given by our theorem is tight for $k=3$ and $k=4$; that is, $b_3=b_{3,5}=6$ and $b_4=b_{4,7}=14$. To this end, let us encode each box $B\subseteq \{0,1\}^n$ as a word $w=w_1w_2\cdots w_n$ over the alphabet $\{0, 1, *\}$. The encoding is defined by the correspondence: $\{0\}\leftrightarrow 0$, $\{1\}\leftrightarrow 1$, $\{0,1\}\leftrightarrow *$. Now, we can describe suitable families of boxes as follows:
$$ 
\begin{array}{rrrrl} 
*& *& 0& 0& 0\\
*& 0& *& 0& 1\\
*& 0& 0& 1& *\\
*& 0& 1& *& 0\\
0& 1& *& *& 1\\
1& 1& *& 1& * \,\,\, ,
\end{array}
$$
$$
\begin{array}{rrrrrrl}
 *& *& 1& 1& 0& 1& *\\
 *& *& 1& 1& 1& *& 0\\
 *& 0& *& 1& 0& 0& *\\
 *& 0& 0& *& 1& 0& *\\
 *& 0& 0& 1& *& 1& *\\
 *& 1& *& 0& 1& *& 0\\
 *& 1& 0& 1& *& *& 0\\
 *& 1& 1& 0& 0& *& *\\
 0& *& 0& 0& 0& *& *\\
 0& 0& 1& 0& *& *& *\\
 0& 1& *& *& 1& *& 1\\
 1& *& 1& *& 1& *& 1\\
 1& 0& *& 0& 0& *& *\\
 1& 1& 0& *& *& *& 1\,\,\,.
\end{array}
$$
It can be easily shown that $b_{k+1}\ge 2b_k$ (see \cite[Sect. 3]{KP}) for all $k$, and thus $b_4=14$ implies $b_k\ge \frac 7 8 2^k$ for $k\ge 4$.  The evaluation (or estimate) of $b_k$ seems to be an interesting problem on its own. The referees guess that $\lim_{k\to \infty}(2^k-b_k)=\infty$.  

\medskip
\noindent
\textbf{2.} Our main result is stronger than \cite[Theorem 1]{KP}. We wish to underline, however, that the methods we have used in order to prove the two results are quite different. At present, we are not able to determine which of them, if any, will be more effective in searching for a more realistic upper bound for the cardinalities of neighbourly families of simplices.

\medskip
\noindent
\textbf{3.} We will not be surprised if Lemma \ref{podstawa} proves to be known. Nevertheless, let us remark that the
 inequality (\ref{nier}) is sharp if $n$ is sufficiently large. It can be shown by the following construction:
\begin{quotation} 
Let $\mathbb N$ be the set of positive integers. We define a sequence $(\mathscr G_v \subset 2^{\mathbb N}\colon v\ge 0)$ by induction. Let $\mathscr G_0=\{ \emptyset, \{1\}\}$. Suppose that $\mathscr G_v$ is already defined. Set $m_v=\max \bigcup \mathscr G_v$. For $A \in \mathscr G'_v=\mathscr G_v\setminus \{\emptyset\}$, let $A^\#= A\cup (A+m_v)$. Let $K=(\bigcup \mathscr G_v+m_v) \cup \{2m_v+1\}$. Then define
$$
\mathscr G_{v+1}=\{A^\#\colon A\in \mathscr G'_v\}\cup \{K\triangle A^\# \colon  A\in \mathscr G'_v\}\cup \{\emptyset, K\}.
$$
(For example, $\mathscr G_1=\{\emptyset, \{1,2\}, \{1,3\}, \{2,3\}\}$.) It is easily checked that $\mathscr  G'_v\subset \binom {\mathbb N}{2^{v}}$, $|\mathscr G_v| =2^{v+1}$ and $m_v=2^{v+1}-1$. Moreover $\ka G_v$ is a group under the symmetric difference. Consequently, for $k=2^v$, $\ka G=\ka G_v$ and $n=2^{v+1}-1$ the inequality (\ref{nier}) becomes an equality. 

For general $k$, we proceed as follows. Write $k$ as a product $k=2^v p$, where $v$ is the $2$-adic order of $k$. (Then $p$ is an odd number.) Define $f\colon \mathbb N\to \mathbb N$ by the formula $f(n)=\lceil \frac n p\rceil$. Set 
$$
\ka G = \{ f^{-1}(A)\colon A \in \ka G_v\}. 
$$
Clearly, it shows that for every $k$, one can find a group $\ka G\subseteq \binom {[n]}{k}$ so that (\ref{nier}) becomes an equality, whenever $n\ge 2k-p$.

\end{quotation}

\end{document}